\documentclass[12pt,a4paper]{article}
\usepackage{amssymb,amsfonts,amsmath,amsthm,amscd,array,hhline}
\DeclareMathOperator*{\Ssym}{\mathrm S}
\DeclareMathOperator*{\rst}{\mathrm rest}
\DeclareMathOperator*{\chrr}{\mathrm char}
\DeclareMathOperator*{\vr}{\mathrm{Var}}

\DeclareMathOperator*{\An}{\mathrm Ann}
\newcommand{\F}{\mathcal F}
\newcommand{\M}{\mathcal V}
\newcommand{\A}{\mathcal A}
\newcommand{\Ag}{\mathfrak A}
\newcommand{\MA}{{\mathcal M}\left(A\right)}
\newcommand{\MAg}{{\mathcal M}\left(\Ag\right)}
\newcommand{\MAgg}{{\mathcal M}\left(\Ag\right)'}
\newcommand{\opL}{L^*}
\newcommand{\tp}[1]{\tau\left(#1\right)}

\newcommand{\Gob}[1]{{\mathrm G}\left(#1\right)}
\newcommand{\Gv}[2]{\Tilde{\mathrm G}_{#2}\left(#1\right)}
\newcommand{\obr}[1]{F\left(#1\right)}
\newcommand{\md}[1]{\left|#1\right|}
\renewcommand{\t}[2]{#1\cdot #2}
\newcommand{\ass}[3]{\left(#1,#2,#3\right)}

\newcommand{\xc}{\Bar{x}}
\newcommand{\crc}[1]{\Bar{#1}}
\newcommand{\smm}[1]{\Check{#1}}

\newcommand{\jord}[2]{#1\circ #2}
\newcommand{\com}[2]{\left[#1,#2\right]}
\newcommand{\dva}[2]{\left(#1#2\right)}
\newcommand{\svob}[1]{\F_{\scriptstyle #1}\left[X\right]}
\newcommand{\svobb}[2]{\F_{\scriptstyle #1}\left[#2\right]}
\newcommand{\svobs}[2]{\F^{(s)}_{\scriptstyle #1}\left[#2\right]}
\newcommand{\pdv}{\mathcal P_d\left(\M\right)}
\newcommand{\pdvv}{\mathcal P_d\left(\mathfrak M\right)}
\newcommand{\ost}[2]{\rst\left(#1,#2\right)}

\newcommand{\lin}[4]{\phi\left(#1,#2,#3,#4\right)}
\newcommand{\linn}[4]{\xi\left(#1,#2,#3,#4\right)}

\newcommand{\cicl}[4]{\psi\left(#1,#2,#3,#4\right)}
\newcommand{\plin}[4]{\vartheta\left(#1,#2,#3,#4\right)}
\newcommand{\linod}[6]{\phi_{#1}\left(#2,#3,#4,#5,#6\right)}
\newcommand{\nlinod}[5]{\phi_{#1}\left(#2,#3,#4,#5\right)}
\newcommand{\plinod}[6]{\vartheta_{#1}\left(#2,#3,#4,#5,#6\right)}
\newcommand{\linodd}[7]{\varphi_{#1}\left(#2,#3,#4,#5,#6,#7\right)}
\newcommand{\nlinodd}[6]{\varphi_{#1}\left(#2,#3,#4,#5,#6\right)}
\newcommand{\chr}[1]{\chrr\left(#1\right)}
\newcommand{\Ann}[1]{\An #1}
\makeatletter
\renewcommand{\section}{\@startsection{section}{1}%
{\parindent}{3.5ex plus 1ex minus .2ex}%
{2.3ex plus.2ex}{\large\bf}}
\renewcommand{\subsection}{\@startsection{subsection}{2}%
{\parindent}{3.25ex plus 1ex minus .2ex}%
{1.5ex plus.2ex}{\normalsize\bf}}
\makeatother
\newtheorem*{theorem*}{Theorem}
\newtheorem{lm}{Lemma}[section]
\newtheorem{pr}{Proposition}[section]
\newtheorem*{rem}{Remark}
\title{\vspace{-24pt}
\bfseries
Non-finitely based varieties of right alternative metabelian algebras}
\author{\bfseries Alexey~Kuz'min\thanks{The author is supported by the
S\~ao Paulo Research Foundation (FAPESP), grant no 2010/51880--2}}
\date{}
\begin{document}
\maketitle

\begin{abstract}
Since 1976, it is known from the paper by V.~P.~Belkin that the variety~$\mathrm{RA_2}$ of right alternative metabelian
(solvable of index~$2$) algebras over an arbitrary field is not Spechtian (contains non-finitely based subvarieties).
In 2005, S.~V.~Pchelintsev proved that the variety generated by the Grassmann
$\mathrm{RA_2}\text{-algebra}$ of finite rank~$r$ over a field~$\mathcal{F}$, for~$\mathrm{char}(\mathcal{F})\neq2$,
is Spechtian iff~$r=1$.
We construct a non-finitely based variety $\mathfrak{M}$ generated by the Grassmann~$\mathcal{V}\text{-algebra}$
of rank~$2$ of certain finitely based subvariety $\mathcal V\subset\mathrm{RA}_2$
over a field~$\mathcal{F}$, for~$\mathrm{char}(\mathcal{F})\neq2,3$,
such that~$\mathfrak{M}$ can also be generated by the Grassmann envelope
of a five-dimensional superalgebra with one-dimensional even part.
\end{abstract}

\begin{list}{}{\rightmargin=\leftmargin\small}
\item\hspace{16pt}{\bfseries Key words:}
non-finitely based variety of algebras,
Spechtian variety of algebras,
right alternative metabelian algebra,
superalgebra,
Grassmann algebra.
\bigskip

\end{list}

\section*{Introduction}
A variety of algebras is said to be \textit{Spechtian}
(or to have the \textit{Specht property}) if its every subvariety is finitely based.
In 1986, A.~R.~Kemer~\cite{Kemer87,Kemer88} solved the famous Specht problem~\cite{Specht50} by proving
that the variety of associative algebras over a field of characteristic~$0$ is Spechtian.
A.~Ya.~Belov~\cite{Belov00},
A.~V.~Grishin~\cite{Grishin99},
and
V.~V.~Schigolev~\cite{Schigolev00}
constructed, independently,
non-finitely based varieties of associative algebras over a field
of prime characteristic.

The Specht property problems for varieties of nonassociative algebras are studied
hard~(see~[1,\,3,\,5--7,\,10--17,\,22--26]).
In 1968, M.~R.~Vaughan-Lee~\cite{Vaughan-Lee68} proved the Specht property of the variety of  metabelian
Lie algebras over a field of characteristic distinct from~$2$.
Also, in his work~\cite{Vaughan-Lee70}, M.~R.~Vaughan-Lee constructed a non-finitely based variety
of Lie algebras over a field of characteristic~$2$.
V.~S.~Drensky~\cite{Drensky74} generalized this result of~\cite{Vaughan-Lee70} to the case of a field of arbitrary prime characteristic.
Yu.~A.~Medvedev~\cite{Medvedev78}
proved the Specht property of the variety of metabelian Malcev algebras.
U.~U.~Umirbaev~\cite{Umirbaev84} generalized this result of~\cite{Medvedev78} to the case of metabelian binary-Lie algebras.
Besides, in his work~\cite{Umirbaev85}, U.~U.~Umirbaev proved the Specht property of every
solvable variety of alternative algebras over a field of characteristic distinct from~$2$ and~$3$.
The essentiality of these restrictions on the characteristic of a ground field is proved by
Yu.~A.~Medvedev~\cite{Medvedev80} and S.~V.~Pchelintsev~\cite{Pchelintsev00}.

There are analogs of the Kemer's Theorem~\cite{Kemer87} in the cases of Jordan, alternative, and Lie algebras
over a field of characteristic~$0$.
Namely, A.~Ya.~Vais and E.~I.~Zel'manov~\cite{Vais-Zelmanov89}
proved that a finitely generated Jordan
$\mathrm{PI}\text{-algebra}$ generates a Spechtian variety.
A.~V.~Iltyakov obtained the similar results for alternative
$\mathrm{PI}\text{-algebras}$~\cite{Iltyakov91} and for finite dimensional Lie algebras~\cite{Iltyakov92}.
Nevertheless, the Specht property problems for the varieties of all alternative, Lie, and Jordan algebras
are still open.

\medskip
Let $\F$ be a field of characteristic distinct from~$2$.
Consider the identities
\begin{align}
\left(x,y,y\right)&=0\quad(\textit{the right alternative identity}),\label{prav}\\
\left(xy\right)\left(zt\right)&=0\quad(\textit{the metabelian identity}),\label{metab}\\
\left(x\circ y\right)\circ z&=0\quad(\textit{the identity of Jordan nilpotency of step 2}),\label{jornil2}
\end{align}
where
$\left(a,b,c\right)=(ab)c-a(bc)$ is the associator of the elements~$a,b,c$
and
$a\circ b=ab+ba$ is the Jordan product of the elements~$a,b$.
The
\textit{variety
$\mathrm{RA_2}$ of right alternative metabelian algebras over $\F$}
is defined by~\eqref{prav},~\eqref{metab}.
By $\mathrm{RA_2'}$ we denote the \textit{subvariety of
$\mathrm{RA_2}$ distinguished by~\eqref{jornil2}}.

Since 1976, it is known~\cite{Belkin76} that
$\mathrm{RA_2}$ is not Spechtian.
I.~M.~Isaev~\cite{Isaev86} proved that non-finitely based subvarieties of
$\mathrm{RA_2}$
can even be generated by finite-dimensional algebras.
Although it was not mentioned by the authors, the direct verification shows that the algebras
constructed in~\cite{Belkin76,Isaev86}
satisfy~\eqref{jornil2}, i.~e. the referred results hold for
$\mathrm{RA_2}'$ as well.
On the other hand, a number of corollaries of the
Yu.~A.~Medvedev's Theorem on two-term identities~\cite{Medvedev78}
states the Specht property of the subvarieties of alternative, left-nilpotent, and
$(-1,1)\text{-algebras}$ in $\mathrm{RA}_2$.
Certain generalizations of these results of~\cite{Medvedev78} are obtained by the author in~\cite{Kuz'min06}.

Recall~\cite{Pchelintsev07} the notion of Grassmann $\M\text{-algebra}$ of finite rank.
Let $\M$ be a variety of algebras over~$\F$;
$\mathcal A=\mathcal A_0\oplus\mathcal A_1$
be a \textit{superalgebra}
($\mathbb Z_2\text{-graded algebra}$) with the
even part $\mathcal A_0$
and the
odd part $\mathcal A_1$,
i.~e.
$\mathcal A_i\mathcal A_j\subseteq\mathcal A_{i+j\;\,(\mathrm{mod} 2)}$
for
$i,j\in\{0,1\}$;
${\mathrm G}$
be the \textit{Grassmann algebra} on a countable set of anticommuting generators
$\left\{e_1,e_2,\ldots\mid e_ie_j=-e_je_i\right\}$
with the natural $\mathbb Z_2\text{-grading}$
($\mathrm G_0$ and $\mathrm G_1$ are spanned by the words of even and, respectively, odd length on the generators
$\left\{e_i\right\}$).
The \textit{Grassmann envelope
${\mathrm G}\left(\mathcal A\right)$
of the superalgebra $\mathcal A$} is the subalgebra
$
{\mathrm G_0}\otimes \mathcal A_0+{\mathrm G_1}\otimes \mathcal A_1
$
of the tensor product
${\mathrm G}\otimes \mathcal A$.
Recall~\cite{Shestakov91,Shestakov93,Vaughan-Lee98,Zelmanov-Shestakov90}
that $\A$ is said to be a
$\M\textit{-superalgebra}$ if ${\mathrm G}\left(\A\right)\in\M$.
Consider a free $\M\text{-superalgebra}$
$\svobs{\M}{U}$ on some finite set
$U=\left\{u_1,\ldots,u_r\right\}$ of free odd generators.
Let us fix in its Grassmann envelope
the elements
$u_{ij}=u_i\otimes e_{i+rj}$,
where
$i=1,\dots,r$ and $j=0,1,\dots$
Then a subalgebra of
${\mathrm G}\bigl(\svobs{\M}{U}\bigr)$
generated by all the elements
$u_{ij}$
is called the \textit{Grassmann $\M\textit{-algebra}$ of rank~$r$}
and is denoted by
$\Gv{\M}{r}$.
We stress that, by definition, the generators of
$\Gv{\M}{r}$
form~$r$ distinct countable families
$\left\{u_{i0},\, u_{i1},\dots, u_{in},\ldots\mid i=1,\dots,r\right\}$
such that every monomial of
$\Gv{\M}{r}$
is skew-symmetric with respect to its variables that belong to the same family.

In 2005, S.~V.~Pchelintsev~\cite{Pchelintsev05,Pchelintsev07} studied the identities of
Grassmann $\mathrm{RA}_2\text{-algebras}$ of ranks~$1$ and~$2$.
In particular, a finite basis for identities of
$\Gv{\mathrm{RA}_2}{1}$ was constructed and
the Specht property of the variety $\vr\Gv{\mathrm{RA}_2}{1}$ generated by
$\Gv{\mathrm{RA}_2}{1}$ was proved.
Moreover, it was shown that
$\vr\Gv{\mathrm{RA}_2}{2}$
is not Spechtian.

In view of the referred results, the following question gives rise:
whether for every finitely based variety
$\M$
its Grassmann algebra
$\Tilde{\mathrm G}_{r}\left(\M\right)$
has a finite basis for identities?
In the present paper, we give the negative answer to this question.

\medskip

Let us fix a field
$\F$
of characteristic
$\chr{\F}\neq2,3$
and consider the
\textit{subvariety $\M$ of $\mathrm{RA_2'}$ distinguished by the identity}
\begin{equation}\label{coscom}
\com{\ass{x}{yz}{x}}{t}=0,
\end{equation}
where
$\com{a}{b}=ab-ba$
is the commutator of the elements~$a,b$.
We prove the following

\smallskip

\begin{theorem*}
The variety
$\vr\Gv{\M}{2}$ is a non-finitely based subvariety of $\M$
distinguished by the system of identities
\begin{equation}\label{kuzass}
\biggl(x,\Bigl(y_1,\ldots,\bigl(y_{n-1},
\ass{y_n}{x}{y_n},y_{n+1}\bigr),\ldots,y_{2n-1}\Bigr),x\biggr)=0,\quad n=1,2,\ldots
\end{equation}
Moreover, $\vr\Gv{\M}{2}$ can be generated by the Grassmann envelope
of a five-dimensional superalgebra on two odd generators
with one-dimensional even part.
\end{theorem*}

\section{Linear generators of the free algebra~$\svob{\M}$}
\label{Sec:AddStr}

\subsection{Common notations}%
Throughout the paper,
all the vector spaces (algebras, superalgebras) are considered over the field~$\F$.
Let us fix the following notations:
\begin{trivlist}
\item[]
$\ost{n}{m}$ is a rest of the integer division of~$n$ by~$m$;
\item[]
$L_a$ and $R_a$ are operators of left and right multiplication by the element~$a$, respectively;
\item[]
$T_a$ is the common denotation for $L_a$ and $R_a$;
\item[]
$\MA$ is an algebra of multiplications of an algebra $A$,
i.~e. an associative algebra that is generated by all the operators $T_a$ $\left(a\in A\right)$ and by the identical mapping
$\mathrm{id}$;
\item[]
$\MA'$ is an algebra generated by the restrictions of all operators from~$\MA$ on $A^2$;
\item[]
$X=\{x_1,x_2,\ldots\}$ is a fixed countable set and
$X_d=\{x_1,x_2,\ldots,x_d\}$;
\item[]
$\svobb{\M}{Y}$ is a free algebra of a variety $\M$ on a set $Y$ of free generators over $\F$;
\item[]
$T_i=T_{x_i}$ is an operator of $\mathcal M\left(\svobb{\M}{X}\right)$;
\item[]
$\pdv$ is a subspace of all multilinear polynomials of degree~$d\geqslant 3$ in $\svobb{\M}{X_d}$;
\item[]
${\mathrm S}_n$ is a symmetric group on the set $1,2,\dots,n$;
\item[]
${\mathrm C}_n$ is a subgroup of ${\mathrm S}_n$ generated by the cycle~$\left(1\,2\ldots n\right)$.
\end{trivlist}
Let
$f=f\left(x_{1},\ldots, x_{n}\right)$
be a polynomial of~$\svob{\M}$
that is linear with respect to some its variables
$x_{i_1},\dots,x_{i_k}$,
$k\geqslant2$.
Then we set
$$
f\left(x_1,\dots,\smm{x}_{i_1},\dots,\smm{x}_{i_k},\dots, x_n\right)=
\sum_{\sigma\in {\mathrm S}_k}
f\bigl(x_1,\dots,x_{i_{\sigma(1)}},\dots, x_{i_{\sigma(k)}},\dots,x_{n}\bigr),
$$
where the permutations are realized with respect to the variables
$x_{i_1},\dots,x_{i_k}$.
The symbol~$\smm{\phantom{x}}$ indicates the variables taking part in the permutations.
Similarly,
$$
f\left(x_1,\dots,\crc{x}_{i_1},\dots,\crc{x}_{i_k},\dots, x_n\right)=
\sum_{\sigma\in {\mathrm C}_k}
f\bigl(x_1,\dots,x_{i_{\sigma(1)}},\dots, x_{i_{\sigma(k)}},\dots,x_{n}\bigr).
$$

While writing down operators of
$\mathcal M\left(\svobb{\M}{X}\right)$
we mark naturals with the symbols~$\Bar{ }$ and~$\Check{ }$
assuming that these symbols are arranged over the variables with the indices equal to the marked naturals.

Let
$\mathcal A=\mathcal A_0\oplus\mathcal A_1$
be a superalgebra.
It is well known that
${\mathrm G}\left(\mathcal A\right)$
satisfies a multilinear identity
$f=0$
iff
$\mathcal A$
satisfies the graded identity
$\Tilde{f}=0$
called the
\textit{superization of~$f=0$.}
Here,
$\Tilde{f}$ denotes the so-called
\textit{superpolynomial corresponding to $f$} and
we say that
$\mathcal A$ satisfies the \textit{superidentity}
$\Tilde{f}=0$.
The detailed descriptions of the process of constructing of superpolynomials
(the \textit{superizing process})
can be found in~\cite{Shestakov91,Shestakov93,Vaughan-Lee98,Zelmanov-Shestakov90}.

\subsection{Operator relations}%
\label{SubSec:SootnosheniyaAlgebry}

\begin{lm}\label{lm jord-bel}
A right alternative algebra $A$ satisfies identity~\eqref{jornil2} if and only if
the operator $L_{a}R_{b}$ is skew-symmetric with respect to $a,b$ in $\MA$.
\end{lm}
\proof
First assume that
$L_{a}R_{b}=-L_{b}R_{a}$ in $\MA$.
Then using~\eqref{prav}, we have
$$
\jord{a^2}{b}=a^2b+ba^2=a^2b+\dva{b}{a}a=
a\bigl(L_{a}R_{b}+L_{b}R_{a}\bigr)=0.
$$
By linearization of the obtained equality we get~\eqref{jornil2}.

Conversely, combining~\eqref{prav} with its linearization
and~\eqref{jornil2},
we obtain
$$
2bL_{a}R_{a}=\dva{a}{b}a-a^2b+a\left(\jord{b}{a}\right)
=\dva{a}{b}a+ba^2-\left(\jord{a}{b}\right)a =ba^2-\dva{b}{a}a=0.\quad\square
$$

\smallskip
We set
$\Ag=\svob{\M}$.
Lemma~\ref{lm jord-bel} yields that $\MAg$ satisfies the relation
\begin{equation}\label{lr}
L_xR_x=0.
\end{equation}

\begin{pr}\label{pr soot zvezd}
The algebra $\MAgg$ satisfies the relations
\begin{align}
&R_xR_y=0\label{rr},\\
&L_xL_y=\com{L_y}{R_x}\label{ll},\\
&R_xL_{x}R_{y}=0\label{RxLxRy}.
\end{align}
\end{pr}

\proof
Suppose
$w\in\Ag^2$.
Using~\eqref{metab} and~\eqref{lr}, we have
$$
wR_xR_y=xL_{w}R_{y}=-xL_{y}R_{w}=-\dva{y}{x}w=0.
$$
Applying the linearization of~\eqref{prav} with~\eqref{metab}, we obtain
$$
wL_xL_y=-\ass{y}{x}{w}=\ass{y}{w}{x}=w\com{L_y}{R_x}.
$$
Combining~\eqref{metab},~\eqref{coscom},~\eqref{lr},~\eqref{rr}, and~\eqref{ll}, we calculate:
\begin{multline*}
2wR_xL_xR_y=2w\com{R_x}{L_x}R_y=-2\ass{x}{w}{x}y
=-\jord{\ass{x}{w}{x}}{y}-\com{\ass{x}{w}{x}}{y}=\\
=\bigl(x\dva{w}{x}\bigr)\left(R_y+L_y\right)
=\dva{w}{x}\left(L_xR_y+L_xL_y\right)
=\dva{w}{x}\left(L_xR_y+\com{L_y}{R_x}\right)=\\
=\dva{w}{x}\left(L_xR_y+L_yR_x-R_xL_y\right)
=-\dva{w}{x}R_xL_y=-wR^2_xL_y=0.\enskip\square
\end{multline*}

\subsection{Standard operators}
\label{SubSec:PravilOper}

A
\textit{standard operator}
is the identical mapping or an operator
$H\in\MAg$ of the form
$$
H=R^{\varepsilon}_{j_0}L_{j_1}R_{j_2}\ldots L_{j_{2k-1}}R_{j_{2k}}L^{\varepsilon'}_{j_{2k+1}},
$$
where
$\varepsilon,\varepsilon'\in\left\{0,1\right\}$
and
$T^{\varepsilon}_{x}=\left\{
\begin{aligned}
\mathrm{id},&\enskip\text{ if}\enskip\varepsilon=0,\\
T_x,&\enskip\text{ if}\enskip\varepsilon=1;
\end{aligned}\right.$
the pair $\left(\varepsilon,\varepsilon'\right)$
is called the
\textit{type of $H$}
and is denoted by
$\tp{H}$.

\bigskip
We stress that applying relations~\eqref{rr} and~\eqref{ll} it is not hard to prove that the following lemma holds.

\begin{lm}\label{lm linagzvezd}
The algebra $\MAg$ is a linear span of standard operators.
\end{lm}

\smallskip
Furthermore, the following lemma is an immediate consequence of relations~\eqref{lr} and~\eqref{RxLxRy}.

\begin{lm}\label{lm cosreg}
A standard operator of type~$\left(0,0\right)$ is skew-symmetric with respect to all its variables.
\end{lm}

\subsection{Standard monomials}
\label{SubSec:PravilnyeOdnochleny}

Let $w$ be a monomial in $\Ag$.
Then $w$ is called the
\textit{standard monomial of type $\tp{w}$}
if
$w=x_{i}H$,
where $H$ is a standard operator distinct from~$R_{j}$
and
$\tp{w}=\tp{H}$.
Note that, by definition,
the elements of
$X^2$ are standard monomials of type~$(0,1)$.
Further,
an
\textit{origin of~$w$}
is a monomial
$w_0=x_iR^{\varepsilon}_{j}$, where
$\tp{w}=\left(\varepsilon,\varepsilon'\right)$;
a
\textit{formative operator of~$w$}
is an operator
$\obr{w}$ such that
$
w=w_0\obr{w}L^{\varepsilon'}_{k}.
$

\bigskip
The following lemma is an immediate consequence of Lemma~\ref{lm linagzvezd}.
\smallskip
\begin{lm}\label{lm linag}
The algebra $\Ag$ is spanned by the standard monomials.
\end{lm}

\medskip

We say that a standard monomial $w$ is \textit{nondegenerate} if $\obr{w}\neq\mathrm{id}$.
Otherwise, $w$ is called \textit{degenerate}.

\smallskip
\begin{lm}\label{lm cosnach}
Every nondegenerate standard monomial~$w$ satisfies the following conditions.
\begin{enumerate}
\item
The formative operator
$\obr{w}$ is skew-symmetric with respect to all its variables.
\item
The origin $w_0$ is skew-symmetric with respect to its variables.
\end{enumerate}
\end{lm}

\proof
The first condition follows from Lemma~\ref{lm cosreg}.
Using~\eqref{prav} and~\eqref{rr}, we prove the second one:
$$
y^2L_{x}R_{z}=\left(xy^2\right)z=\left(xy\right)R_yR_z=0.\quad\square
$$

\subsection{Basis monomials}
\label{SubSec:BazisOdn}

A  \textit{basis monomial}
is a standard monomial $w$ such that
the sequences of indices of the variables of the origin~$w_0$
and of the formative operator~$\obr{w}$ ascend strictly.
In particular,
all the standard monomials of degrees~$1$ and~$2$ are basis ones.

\smallskip
\begin{lm}\label{lm porogd odnochleny}
The algebra $\Ag$ is spanned by the basis monomials.
\end{lm}

\proof
By virtue of lemma~\ref{lm linag}, it suffices to prove that
every standard monomial of degree not less then~$3$
can be represented as a linear combination of basis monomials.
Let us rewrite down~\eqref{prav} in the form
$
y^2L_x=yL_xR_y.
$
Then it is clear that an origin of a degenerate standard monomial of degree~$3$
is skew-symmetric modulo linear combinations of nondegenerate standard monomials.
Hence, to conclude the proof it remains to note that by lemma~\ref{lm cosnach},
every nondegenerate standard monomial is proportional to a basis one.
\endproof

\medskip

A \textit{basis polynomial}
is a linear combination over~$\F$ of pairwise distinct basis monomials with nonzero scalars.

\medskip

\begin{lm}\label{lm TAgMultlin}
Every $\mathrm{T}\text{-ideal}$ of $\Ag$ can be generated by
a system of multilinear basis polynomials.
\end{lm}

\proof
Note that every basis monomial by its definition
has a degree not more then~$3$ with respect to any of its variables.
Consequently, in view of the restrictions
$\chr{\F}\neq2,3$,
a
$\mathrm{T}\text{-ideal}$ of $\Ag$
generated by some basis polynomial
can be also generated by a system of multilinear polynomials~(see~\cite[Chap.~1]{Zhevlakov-Slin'ko-Shestakov-Shirshov}).
Moreover, by Lemma~\ref{lm porogd odnochleny}, this polynomials can be expressed linearly with
multilinear basis polynomials.
\endproof

\section{Auxiliary $\M\text{-superalgebra}$}%
\label{Sec:VspomogatelnyeSuperalgebry}

Let $\A=\A_0\oplus\A_1$ be a superalgebra
$$
\A_0=\F\cdot a,\quad\;
\A_1=\F\cdot v+\F\cdot w+\F\cdot y+\F\cdot z+\F\cdot z'
$$
such that all nonzero products of its basis elements are the following:
$$
\t{z}{z}=a,\quad\t{y}{a}=v,\quad\t{z'}{a}=w,\quad\t{y}{v}=\t{v}{y}=a.
$$
By definition,
it is not hard to see that
$\A$ is a metabelian superalgebra generated by the odd elements~$y,z,z'$
and
$w$ lies in the annihilator of $\A$.
Moreover, the direct verification shows that the following proposition holds.

\begin{pr}\label{pr Sup Rel}
The superalgebra
$\mathcal A$
satisfies the relations
\begin{align}
\A_0\A&=0,\label{usl1}\\
\com{\A_1}{\A_1}&=0,\label{usl2}\\
\A_0\com{L_{u_i}}{L_{u_j}}&=0,\label{usl3}\\
\A_1\A_0\com{L_{u_i}}{L_{u_j}}L_{u_k}&=0,\label{usl4}
\end{align}
where
$u_i,u_j,u_k$ are arbitrary generators of $\A$.
\end{pr}

\medskip

We stress that relation~\eqref{usl1} yields the following

\begin{lm}\label{lm LR anih A1}
The operator
$L_{a}R_{b}$
annihilates $\A_1$ for all~$a,b\in\A$.
\end{lm}

\smallskip

\begin{lm}\label{lm dostMsuper}
Every metabelian superalgebra
$\mathcal A$
generated by odd elements
and satisfying~\eqref{usl1}--\eqref{usl4}
is a $\M\text{-superalgebra}$.
\end{lm}

\proof
Throughout the proof,
the conditions of the metability of $\A$,
the odd parity of all generators of~$\A$,
and its consequence
$\A_0\subset\A^2,$
are used with no comments.

First, let us prove that
$\Gob{\mathcal A}$
satisfies~\eqref{prav},
i.~e.
$$
\ass{a}{b}{c}+{\left(-1\right)}^{\md{b}\md{c}}\ass{a}{c}{b}=0,
$$
for all basis elements
$a,b,c\in\A$,
where
$\md{a}$ denotes the parity of the element~$a$
($\md{a}=i\in\left\{0,1\right\}$ if $a\in\A_i$).
Note that if
$\ass{a}{b}{c}\neq0$,
then at least two of the elements
$a,b,c$ are odd.
If $b,c\in\mathcal A_1$, then by~\eqref{usl1},~\eqref{usl2},
the associator
$\ass{a}{b}{c}=a\left(bc\right)$
is symmetric with respect to~$b,c$.
Hence, it remains to check the skew-symmetry of
$\ass{a}{b}{c}$ with respect to $b,c$ under the conditions:
$b\in\A_0$ and $a,c$ are generators of $\A$.
In this case,
using~\eqref{usl1},~\eqref{usl2}, and~\eqref{usl3}, we obtain
$$
\ass{u_i}{b}{u_j}=\left(u_i b\right)u_j=u_j\left(u_ib\right)
=u_i\left(u_jb\right)=-\ass{u_i}{u_j}{b}.
$$
Thus,
$\Gob{\mathcal A}$
is right alternative.

Now, let us prove that
$\Gob{\A}$
satisfies~\eqref{jornil2}.
By Lemma~\ref{lm jord-bel}, it suffices to verify that the operator
$L_{a}R_{b}$ is skew-symmetric in
${\mathcal M}\left(\Gob{\A}\right)$.
Taking into account Lemma~\ref{lm LR anih A1},
it remains to check that
$$
\A_0\left(L_{u_i}R_{u_j}-L_{u_j}R_{u_i}\right)=0.
$$
Indeed, assuming $b\in\A_0$ and applying~\eqref{usl2},~\eqref{usl3}, we have
$$
\left(u_i b\right)u_j=u_j\left(u_ib\right)
=u_i\left(u_jb\right)=\left(u_jb\right)u_i.
$$

Finally, let us prove that
$\Gob{\A}$
satisfies~\eqref{coscom}.
In view of~\eqref{prav} and~\eqref{metab} it suffices to verify that
$L_xL_xR_t=L_xL_xL_t$ in
${\mathcal M}\left(\Gob{\A}\right)'$, i.~e.
$$
b\com{L_{u_i}}{L_{u_j}}R_{u_k}={(-1)}^{\md{b}}b\com{L_{u_i}}{L_{u_j}}L_{u_k},\quad b\in\A^2.
$$
If $b\in\A_0$, then the left side of the equality is zero by virtue of~\eqref{usl1}
and the right side is zero in view of~\eqref{usl3}.
Otherwise, if
$b\in\A_1\cap\A^2=\A_1\A_0$,
then the both sides of the equality are
zeros by virtue of~\eqref{usl2} and~\eqref{usl4}.
\endproof

\smallskip
\begin{lm}\label{lm znachreg}
A value taken by an arbitrary nondegenerate basis monomial $w$ of type~$(\varepsilon,\varepsilon')$ on basis elements of~$\A$
lies in the homogeneous component~$\A_{\varepsilon'}$.
\end{lm}

\proof
Let $0\neq\Tilde{w}\in\A$
be a monomial obtained by a substitution of the variables of~$w$
by arbitrary basis elements of~$\A$.
Then Lemma~\ref{lm LR anih A1} implies
$\Tilde{w_0}\in\A_0$.
Hence, by virtue of the metability and the odd parity of the generators of~$\A$,
we have
$\Tilde{w}_0\obr{\Tilde{w}}\in\A_0$.
Therefore, $\Tilde{w}\in\A_{\varepsilon'}$.
\endproof

\section{Additive basis of the space~$\pdv$}%
\label{Sec:AddBasPdv}

A \textit{regular polynomial}
is a basis polynomial~$f\in\pdv$
represented as a linear combination of nondegenerate basis monomials of a same fixed type.
This type is called the
\textit{type of~$f$} and
is denoted by~$\tp{f}$.
If
$f=0$
is an identity of some algebra
$A\in\M$,
then we say that $A$ satisfies a
\textit{regular identity of type~$\tp{f}$.}

\subsection{Reduction to the regular identities of $\Gob{\A}$}%
\label{SubSec:RedukKRegINereg}

Let
$\mathcal A$
be the
$\M\text{-superalgebra}$ defined in Sec.~\ref{Sec:VspomogatelnyeSuperalgebry}.

\smallskip
\begin{lm}\label{lm raspad}
If
$\Gob{\A}$
satisfies a nontrivial in $\M$ multilinear identity of degree~$d\geqslant3$,
then $\Gob{\A}$ satisfies some regular identity.
\end{lm}

\proof
By virtue of Lemma~\ref{lm porogd odnochleny} we may assume that
$\Gob{\A}$
satisfies an identity
$f=0$,
where
$f\in\pdv$
is a linear combination of pairwise distinct basis monomials with nonzero scalars.

First, consider the case $d=3$.
Suppose that $f$ is not regular of type~$(0,0)$,
i.~e. $f$ contains degenerate basis monomials.
Then, in view of~\eqref{rr},
$fR_{4}$
will be regular of type~$(1,0)$.

In the case $d\geqslant4$, if $f$ is not regular, then it can be represented as a sum
$f_0+f_1$
of two regular polynomials of types
$\tp{f_0}=\bigl(1-\varepsilon,0\bigr)$
and
$\tp{f_1}=\bigl(\varepsilon,1\bigr)$,
where
$\varepsilon=\ost{d}{2}$.
By Lemma~\ref{lm znachreg}, all the values taken by~$\Tilde{f_i}$ $(0,1)$
on basis elements of~$\A$ lie in~$\A_{i}$.
This yields that
$\Gob{\A}$
satisfies both identities~$f_i=0$.
\endproof

\subsection{Reduction to the regular identities of $\Gob{\A}$ of type~$\left(\varepsilon,0\right)$}%
\label{SubSec:RedukciyaKRegulyarnomuSluchayu2}

For
$k=1,\ldots,d+1$
we define linear mappings
$$
\opL_{k}:{\mathcal P}_{d}(\M)\mapsto{\mathcal P}_{d+1}(\M)
$$
acting on the monomials
$w=w(x_1,\ldots,x_{d})\in{\mathcal P}_{d}(\M)$
as follows:
$$
w\opL_{k}=w(x_{1^{\left<k\right>}},\ldots,x_{d^{\left<k\right>}})L_{k},\quad
i^{\left<k\right>}=\left\{
\begin{aligned}
i,\enskip i<k,\\
i+1,\enskip i\geqslant k.
\end{aligned}\right.
$$

\smallskip
\begin{lm}\label{lm raspad0}
If
$\Gob{\A}$
satisfies a regular identity of type~$(\varepsilon,1)$,
then $\Gob{\A}$ satisfies a regular identity of type~$(\varepsilon,0)$.
\end{lm}

\proof
Let $g$ be a regular polynomial of
${\mathcal P}_{d+1}(\M)$
of type~$(\varepsilon,1)$.
Then by the definition of standard monomials we can represent
$g$
in the form
$$
g=\sum_{k\in I}f_k\opL_{k},
$$
where
$\varnothing\neq I\subseteq\{1,\dots,d+1\}$
and every
$f_k\in{\mathcal P}_{d}(\M)$ is a regular polynomial of type~$(\varepsilon,0)$.

Let us prove that the identity $g=0$ in~$\Gob{\A}$ implies~$f_k=0$.
Indeed, assume that
$\Tilde{g}=0$ in~$\A$
and
$\Tilde{f_k}$
takes a nonzero value at some basis elements
$x_1=b_1,\ldots,\;x_{d}=b_{d}$ of~$\A$.
Then Lemma~\ref{lm znachreg} implies~$\Tilde{f_k}\in\A_{0}$
and, consequently, by the definition of~$\A$,
we have
$\Tilde{f_k}=\alpha\,a$,
where
$0\neq\alpha\in\F$.
But in this case,
$\Tilde{g}$
takes a nonzero value in~$\A$ at the elements
$
x_{1^{\left<k\right>}}=b_1,\ldots,
\; x_{d^{\left<k\right>}}=b_{d},\; x_{k}=z':
$
$$
\Tilde{g}=
\Tilde{f_k}L_{z'}=
\alpha\,\t{z'}{a}=\alpha\,w\neq0.
$$
The obtained contradiction completes the proof.
\endproof

\subsection{Linear independence of the basic monomials in $\pdv$}%
\label{SubSec:RegulyarnyiSluchai}

In order to avoid complicated formulas
while writing down the elements of
$\pdv$
we omit the indices of variables at the operator symbols~$L,R$ and assume that they are arranged in the ascending order.
For example, the notation
$\left(x_2x_5\right)\left(LR\right)^2$
means  the monomial
$\left(x_2x_5\right)L_{1}R_{3}L_{4}R_{6}$.

\smallskip

\begin{lm}\label{lm AddBasPd}
The set of all basis monomials in $\pdv$ forms its additive basis.
\end{lm}

\proof
By virtue of Lemmas~\ref{lm porogd odnochleny},~\ref{lm dostMsuper},~\ref{lm raspad} and~\ref{lm raspad0},
it suffices to prove that
$\Gob{\A}$
does not satisfy any regular identity of type~$(\varepsilon,0)$.
Consider a regular polynomial $f\in\pdv$ of type~$\tp{f}=(\varepsilon,0)$.
If $d$ is even, then we can represent $f$ in the form
$$
f=\sum_{i<j}
\alpha_{i,j}\,(x_{i}x_{j}){\left(LR\right)}^{\frac{d}{2}-1},\quad 0\neq\alpha_{i,j}\in\F.
$$
By
$\Tilde{f}_{i,j}$
denote the value taken by the superpolynomial
$\Tilde{f}$
on the following elements of~$\A$:
$$
x_i=x_j=z,\quad x_k=y\;\text{ for  all }\; k\neq i,j.
$$
Hence, $\Tilde{f}_{i,j}$ is proportional to the element
$$
(\t{z}{z}){\left(L_{y}R_{y}\right)}^{\frac{d}{2}-1}=a\neq0.
$$
For odd $d$, we have
$$
f=\sum_{i}
\alpha_{i}\,x_{i}{\left(LR\right)}^{\frac{d-1}{2}},\quad 0\neq\alpha_i\in\F.
$$
Similarly, a value
$\Tilde{f}_i$
taken by
$\Tilde{f}$
on the elements
$$
x_i=a,\quad x_j=y\;\text{ for  all }\; j\neq i
$$
turns out to be proportional to
$$
a{\left(L_{y}R_{y}\right)}^{\frac{d-1}{2}}=a\neq0.\quad\square
$$

\bigskip

We stress that, in view of Lemma~\ref{lm TAgMultlin}, the proof of Lemma~\ref{lm AddBasPd} implies
$\M=\vr\Gv{\M}{3}$.

\section{Auxiliary polynomials}%
\label{Sec:VspomogatelnyeFunkcii}

\subsection{Polynomials $\xi,\psi,\phi$}%
\label{SubSec:FunkciiPsiPhiXi}

Consider the following polynomials in~$\Ag$:
\begin{gather*}
\linn{x}{y}{z}{t}=
\dva{x}{y}L_{z}R_{t}+\dva{z}{t}L_{x}R_{y},\\
\cicl{x}{y}{z}{t}=
\left(\xc\crc{y}\right)L_{\crc{z}}R_{t},\qquad
\lin{x}{y}{z}{t}=
\left(\xc\crc{y}\right)L_{\crc{z}}R_{\crc{t}}.
\end{gather*}
Lemma~\ref{lm cosnach} yields immediately the following properties of~$\xi$ and~$\psi$.
\begin{lm}\label{lm SvXi}
The polynomial
$\linn{x}{y}{z}{t}$
is skew-symmetric w.r.t. the pair~$x,y$ and, independently, w.r.t. the pair~$z,t$.
\end{lm}

\medskip

\begin{lm}\label{lm SvPsi}
The polynomial
$\cicl{x}{y}{z}{t}$
is skew-symmetric w.r.t.
$x,y,z$.
\end{lm}

\medskip

Moreover, combining the definition of~$\xi$ with Lemma~\ref{lm SvXi}, we obtain the following

\begin{lm}\label{lm KleinFour}
The polynomial~$\xi$ is invariant under the action of the Klein four-group on its variables:
$$
\linn{x}{y}{z}{t}=\linn{y}{x}{t}{z}=\linn{t}{z}{y}{x}.
$$
\end{lm}

\medskip

\begin{pr}\label{pr IdentPhiXi}
The algebra $\Ag$ satisfies the identities
\begin{align}
&\lin{ab}{x}{y}{z}=0,\label{eq phi-kv}\\
&\lin{x}{y}{z}{t}
=\linn{\crc{x}}{y}{\crc{z}}{t}=\linn{x}{\crc{y}}{z}{\crc{t}},\label{eq phi-xi}\\
&\cicl{a}{x}{b}{x}+\frac{1}{2}\lin{a}{x}{b}{x}=
2\dva{a}{x}L_{b}R_{x},\label{eq cos-psi}\\
&\lin{a}{\crc{x}}{\crc{y}}{\crc{z}}=0.\label{eq in-cyc}
\end{align}
\end{pr}
\proof
Applying~\eqref{metab},~\eqref{rr},~\eqref{ll}, and Lemma~\ref{lm cosreg}, we have
$$
\lin{ab}{x}{y}{z}
=(ab)\left(R_{x}L_{y}R_{z}+L_{z}L_{x}R_{y}\right)
=(ab)\left(R_{x}L_{y}R_{z}-R_{z}L_{x}R_{y}\right)=0.
$$

To prove~\eqref{eq phi-xi}, first note that the equality
$\linn{\crc{x}}{y}{\crc{z}}{t}=\linn{x}{\crc{y}}{z}{\crc{t}}$
follows from the definition of~$\xi$.
Then we stress that~$\phi$ can be represented, by definition, as follows:
$$
\lin{x}{y}{z}{t}
=\linn{x}{y}{z}{t}
+\linn{t}{x}{y}{z}.
$$
Therefore, applying Lemma~\ref{lm KleinFour} to the second summand, we obtain~\eqref{eq phi-xi}.

Further, combining Lemma~\ref{lm cosnach} with~\eqref{eq phi-xi}, we prove~\eqref{eq cos-psi}:
$$
\cicl{a}{x}{b}{x}+\frac{1}{2}\lin{a}{x}{b}{x}
=\dva{a}{x}L_{b}R_{x}+\dva{x}{b}L_{a}R_{x}+\linn{a}{x}{b}{x}=2\dva{a}{x}L_{b}R_{x}.
$$

Finally, using~\eqref{eq phi-xi} and taking into account Lemma~\ref{lm SvXi}, we get
$$
\lin{a}{\crc{x}}{\crc{y}}{\crc{z}}
=\linn{a}{\crc{x}}{y}{\crc{z}}+
\linn{a}{\crc{y}}{z}{\crc{x}}+
\linn{a}{\crc{z}}{x}{\crc{y}}
=\linn{a}{\smm{x}}{\smm{y}}{\smm{z}}=0.\quad\square
$$

\medskip

The  following lemma is an immediate consequence of~\eqref{eq phi-xi}.
\begin{lm}\label{lm SvPhi}
The polynomial
$\lin{x}{y}{z}{t}$
is symmetric w.r.t. the pair~$x,z$ and, independently, w.r.t. the pair $y,t$.
\end{lm}

\subsection{Special regular polynomials}%
\label{SubSec:SpecialWords}

By definition, the space $\pdv$ can be divided into two components
$$
{\mathcal P}_{d}(\M)={\mathcal P}^{(0,1-\varepsilon)}_{d}(\M)+{\mathcal P}^{(1,\varepsilon)}_{d}(\M),
\quad \varepsilon=\ost{d}{2},
$$
where
${\mathcal P}^{(\varepsilon,\varepsilon')}_{d}$
denotes the subspace of all regular polynomials in $\pdv$ of type~$(\varepsilon,\varepsilon')$.

Let us define some special polynomials in~${\mathcal P}^{(1,1)}_{d}(\M)$
for $d=2n+3$,
$n\in\mathbb N$.

We use $\vartheta$ as a common denotation for the symbols~$\xi,\psi,\phi$.
The polynomial
$\plin{x_i}{x_j}{x_k}{x_{\ell}}$
is denoted shortly by
$
\plin{i}{j}{k}{\ell}.
$

A
\textit{$\vartheta\text{--word}$  of order~$n$}
is a polynomial
$f\in{\mathcal P}^{(1,1)}_{d}(\M)$
of the form
$$
f=\plin{i}{j}{k}{\ell}{(LR)^{n-1}}L_{m}
$$
denoted by
$
f=\plinod{n}{i}{j}{k}{\ell}{m}.
$
The polynomial
$f_0=\plin{i}{j}{k}{\ell}$
is called the \textit{origin of~$f$.}

A
\textit{double $\phi\text{--word}$ of order~$n\geqslant2$}
is a polynomial
$f\in{\mathcal P}^{(1,1)}_{d}(\M)$ $(d\geqslant7)$
of the form
$$
f=\lin{i}{j}{k}{\ell}L_{\Bar{m}}R{(LR)^{n-2}}L_{\Bar{q}}
$$
denoted by
$f=\linodd{n}{i}{j}{k}{\ell}{m}{q}$.
The polynomial
$f_0=\lin{i}{j}{k}{\ell}$
is called the
\textit{origin of~$f$}.
Note that in view of Lemma~\ref{lm cosreg},
$f$
can be represented as a linear combination of two
$\phi\text{--words}$ of order~$n$
with the same origins~$f_0$.

A
\textit{triple $\phi\text{--word}$ of order~$n$}
is a polynomial in~${\mathcal P}^{(1,1)}_{d}(\M)$
of the form
$\linod{n}{\Bar{\imath}}{j}{\Bar{k}}{\ell}{\Bar{m}}$.

\medskip

\begin{pr}\label{pr Ident Doub-Trip Phi}
The following identity holds for $n\geqslant2$:
\begin{equation}\label{eq Doub-Trip Phi}
\linodd{n}{1}{2}{3}{\Bar{4}}{\Bar{5}}{\Bar{6}}=-\linod{n}{\Bar{4}}{1}{\Bar{5}}{3}{\Bar{6}}.
\end{equation}
\end{pr}

\proof
By Lemma~\eqref{lm cosnach}, we have
\begin{multline*}
\lin{a}{b}{c}{x}L_xR_y=
(xa)L_bR_cL_xR_y+(cx)L_aR_bL_xR_y=\\
=-(xa)L_xR_cL_bR_y-(xc)L_xR_aL_bR_y
=-\linn{x}{a}{x}{c}L_bR_y.
\end{multline*}
Multiplying the both sides of the obtained equality by
$HL_x$,
where
$$
H=\left\{
\begin{aligned}
&\mathrm{id},\enskip\text{ if }\;n=2,\\
&L_{8}R_{9}\ldots L_{2n+2}R_{2n+3},\enskip\text{ if }\;n\geqslant 3,
\end{aligned}\right.
$$
we get
$$
\lin{a}{b}{c}{x}L_xR_yHL_x=-\linn{x}{a}{x}{c}L_bR_yHL_x.
$$
Hence, by the linearization
$$
a=x_1,\quad b=x_2,\quad c=x_3,\quad
x\mapsto x_4,x_5,x_6,\quad y=x_7,
$$
taking into account~$\eqref{eq phi-xi}$, we obtain
$$
\linodd{n}{1}{2}{3}{\Bar{4}}{\Bar{5}}{\Bar{6}}
=-\xi_{n}\left(\Check{4},1,\Check{5},3,\Check{6}\right)
=-\linod{n}{\Bar{4}}{1}{\Bar{5}}{3}{\Bar{6}}.\quad\square
$$

\bigskip

Let us set
$$
g_n=g_n\left(x,y_1,\dots,y_{2n-1}\right)=
\biggl(x,\Bigl(y_1,\ldots,\bigl(y_{n-1},
\ass{y_n}{x}{y_n},y_{n+1}\bigr),\ldots,y_{2n-1}\Bigr),x\biggr).
$$

\smallskip

\begin{lm}\label{lm Lin gn}
A linearization of~$g_n$ in
${\mathcal P}^{(1,1)}_{d}(\M)$
is proportional to a triple
$\phi\text{--word}$ of order~$n$.
\end{lm}

\proof
Applying Lemma~\ref{lm cosreg} and relation~\eqref{rr}, we have
\begin{multline*}
g_n\left(x,z_{2n-2},\dots,z_2,y,z_1,z_3,\dots,z_{2n-3}\right)=\\
=x\left[L_y,R_y\right]\left[L_{z_2},R_{z_1}\right]\dots\left[L_{z_{2n-2}},R_{z_{2n-3}}\right]\left[L_x,R_x\right]=\\
\shoveright{=xR_yL_y\left[L_{z_2},R_{z_1}\right]\dots\left[L_{z_{2n-2}},R_{z_{2n-3}}\right]R_xL_x=}\\
\shoveleft{={(-1)}^{n-1}(xy)L_yR_{z_1}L_{z_2}\dots R_{z_{2n-3}}L_{z_{2n-2}}R_xL_x=}\\
={(-1)}^{n}(xy)L_xR_yL_{z_1}R_{z_2}\dots L_{z_{2n-3}}R_{z_{2n-2}}L_x.
\end{multline*}
Hence, linearizing
$g_n\mapsto\Delta g_n\in {\mathcal P}^{(1,1)}_{d}(\M)$:
$$
x\mapsto x_1,x_2,x_3,\quad
y\mapsto x_4,x_5,\quad
z_1=x_6,\dots,\;z_{2n-2}=x_{d},
$$
and taking into account~\eqref{eq phi-xi}, we obtain
\begin{equation}\label{eq Lin gn}
\Delta g_n={(-1)}^{n}\xi_n\left(\Check 1,4,\Check 2,5,\Check 3\right)
={(-1)}^{n}\linod{n}{\Bar{1}}{4}{\Bar{2}}{5}{\Bar{3}}.\quad\square
\end{equation}

\begin{lm}\label{lm IntersecT}
The intersection
$
{\mathcal I} ={\left(g_n\right)}^{\mathrm T}\cap
\Bigl(\bigcup_{d=2n+3}^{\infty}\pdv\Bigr)
$
is spanned by the linearizations of~$g_n$.
\end{lm}

\proof
Taking into account Lemma~\ref{lm Lin gn}, it remains to prove that
$$
{\left(\Delta g_n\right)}^{\mathrm T}\cap
\Bigl(\bigcup_{d=2n+4}^{\infty}\pdv\Bigr)=\{0\}.
$$
In view of~\eqref{metab} and~\eqref{eq phi-kv}, it suffices to verify that
a triple
$\phi\text{--word}$ of any order
lies in the annihilator $\Ann \Ag$.
Following the proof of Lemma~\ref{lm Lin gn},
we may restrict with checking
$h\in\Ann \Ag$
for the monomial
$$
h=(xy)L_xR_y{(LR)}^{n-1}L_x.
$$
Indeed,
Lemma~\ref{lm cosreg} yields immediately
$hR_z=0$
and,
using~\eqref{rr} and~\eqref{ll},
we get
$$
hL_z=(xy)L_xR_y{(LR)}^{n-1}L_xL_z=(xy)L_xR_y{(LR)}^{n-1}L_zR_x=0.\quad\square
$$

\section{Linear generators of the space~$\pdvv$}
\label{Sec:AddStrPdvv}

Let $\mathfrak M$ be a subvariety of $\M$ distinguished  by system~\eqref{kuzass}.
In what follows, considering the free algebra
$\mathfrak A'=\svob{\mathfrak M}$ and
its subspace
$\pdvv$,
we assume that they inherit naturally all the notions introduced for~$\Ag$ and~$\pdv$ in previous sections.


\subsection{Normalized words}%
\label{SubSec:NormalWords}

Let $\varPhi_d$ be a
\textit{linear span of all
$\phi\text{--words}$ in~${\mathcal P}^{(1,1)}_{d}(\mathfrak M)$}
for $d=2n+3$,
$n\in\mathbb N$.
A
$\phi\text{--word}$
$f\in\varPhi_d$ of the form
$$
f=\lin{i}{j}{k}{\ell}{(LR)^{n-1}}L_{m}
$$
is called
\textit{normalized}
if
$m$ is a maximal index in the operator
${(LR)^{n-1}}L_{m}$.
In particular, every $\phi\text{--word}$ of order~$1$ is normalized.
We denote a normalized
$\phi\text{--word}$
shortly,
omitting the corresponding minimal index:
$$
f=\nlinod{n}{i}{j}{k}{\ell}.
$$
In view of Lemma~\ref{lm cosreg}, the given definition implies instantly the following

\smallskip

\begin{lm}\label{lm PhiSpanNormal}
Every
$\phi\text{--word}$
$f$
is either normalized or can be expressed linearly with a normalized
$\phi\text{--word}$
and a double
$\phi\text{--word}$
with the same origins~$f_0$.
\end{lm}

\medskip

Let
$\varPhi'_{d}$, for $d\geqslant7$,
be a
\textit{subspace of $\varPhi_{d}$ generated by all double
$\phi\text{--words}$.}
A  double $\phi\text{--word}$
$f\in\varPhi'_d$
of the form
$$
f=\lin{i}{j}{k}{\ell}L_{\Bar m}{R(LR)^{n-2}}L_{\Bar q},\quad n\geqslant2
$$
is called
\textit{normalized}
if
$m$ is a minimal index in the operator
$L_{m}R{(LR)^{n-1}}$.
We denote
$f$
shortly,
omitting $m$:
$$
f=\nlinodd{n}{i}{j}{k}{\ell}{q}.
$$

\smallskip

\begin{lm}\label{lm DoublePhiSpanNormal}
Every double
$\phi\text{--word}$
$f=\linodd{n}{i}{j}{k}{\ell}{m}{q}$
is either normalized or can be expressed linearly with two normalized double
$\phi\text{--words}$ with the same origins~$f_0$.
\end{lm}
\proof
Let $m'$ be a minimum of the set
$\{1,\dots,d\}\setminus \{i,j,k,\ell\}$.
By assumption of the lemma it is clear that
$m'<m$ and $m'<q$.
We stress that in view of Lemma~\ref{lm cosreg},
$f$ satisfies the assertion of Lemma~\ref{lm J-triples}
for~$x_m,x_q,x_{m'}$.
Consequently, we can represent $f$ in the form
$$
f=
\pm\nlinodd{n}{i}{j}{k}{\ell}{m} \pm \nlinodd{n}{i}{j}{k}{\ell}{q}.\quad\square
$$

We call the procedure described in Lemma~\ref{lm DoublePhiSpanNormal} the
\textit{normalization of double
$\phi\text{--word}$.}

\subsection{Tame words}%
\label{SubSec:TameWords}

Let
$
f=f(x,y,\dots,z)
$
be a nonassociative multilinear polynomial and
$
\Ssym{\{f\}}
$
be the symmetric group on the set
$
x,y,\dots,z
$.
For
$
\sigma\in \Ssym{\{f\}}
$
we set
$$
f^{\sigma}=f\bigl(\sigma(x),\sigma(y),\dots,\sigma(z)\bigr).
$$
The key part in the proof of the statements of this subsection is played by the following obvious lemma.

\smallskip

\begin{lm}\label{lm J-triples}
If
$f=f\left(x,y,z_1,\dots,z_n\right)$
is symmetric w.r.t. $x,y$ and
$$
f\left(\crc{x},\crc{y},\crc{z}_1,z_2,\dots,z_n\right)=
f\left(\crc{x},\crc{y},z_1,\crc{z}_2,\dots,z_n\right)=
f\left(\crc{x},\crc{y},z_1,z_2,\dots,\crc{z}_n\right)=0,
$$
then the vector space
$\mathrm{Vec}_{\F}\left<f^{\sigma}\mid \sigma\in \Ssym{\{f\}}\right>$
is spanned by the elements
$f^{\sigma}$
such that
$\sigma(x)=x$.
\end{lm}

\medskip

A $\phi\text{--word}$
$f\in\varPhi_d$
is called
\textit{tame}
if its origin $f_0$ has one of the following types:
$$
1)\enskip\lin{1}{i}{j}{d}, \enskip\, i<j;\qquad\;
2)\enskip\lin{1}{2}{d}{j}.
$$
Otherwise,
$f$ is called
\textit{wild}.

\begin{lm}\label{lm PhiSpanTame}
The space $\varPhi_d$ is spanned by the tame
$\phi\text{--words}$.
\end{lm}

\proof
Let us show that every wild $\phi\text{--word}$
$f\in\varPhi_d$
can be expressed linearly with tame $\phi\text{--words}$.
We do it at four steps, referring each time Lemma~\ref{lm J-triples}
and using
Lemmas~\ref{lm cosreg},~\ref{lm SvPhi}
with no comments.

First let us show that $f$ is a linear combination of $\phi\text{--words}$ with origins
of the form
$
\lin{1}{i}{j}{k}.
$
Indeed, as far as $\phi$ is cyclic, it suffices to assume that $f_0$ doesn't contain the variable~$x_1$.
In this case, taking into account that
every triple $\phi\text{--word}$, in view of Lemma~\ref{lm IntersecT}, is zero in $\pdvv$,
we see that the assertion of Lemma~\ref{lm J-triples} holds
for the variables on the first and the third positions in~$f_0$ and for all the variables outside $f_0$.
Thus applying Lemma~\ref{lm J-triples} under the assumption $x=x_1$, we represent $f$ in the required form.

Further, by the similar arguments for the variables on the second and the fourth positions in~$f_0$,
one can show that if $f_0$ doesn't contain $x_d$, then
$f$ is a linear combination of $\phi\text{--words}$ with origins
of the form
$
\lin{1}{i}{j}{d}.
$

Now suppose that $f$ is wild and
$
f_0=\lin{1}{i}{d}{j}.
$
Then applying Lemma~\ref{lm J-triples} under the assumption $x=x_2$, we can express
$f$ with two tame $\phi\text{--words}$ with the origins of type~2).

Finally, consider the case
$
f_0=\lin{1}{j}{i}{d},
$
where $j>i$.
Using~\eqref{eq in-cyc} and, if necessary, Lemma~\ref{lm J-triples},
we express $f$ with one $\phi\text{--word}$ with the origin of type~1)
and one (in the case $i=2$) or two $\phi\text{--words}$ with the origins of type~2).
\endproof

\medskip

A double $\phi\text{--word}$ in
$\varPhi'_{d}$ is called
\textit{tame}
if it is normalized and has one of the  following types:
$$
1)\enskip\nlinodd{n}{1}{i}{j}{d}{k}, \enskip\,i<j<k;\qquad
2)\enskip\nlinodd{n}{1}{2}{d}{j}{k}, \enskip\,j<k.
$$
Otherwise,
$f$ is called
\textit{wild}.

\smallskip

\begin{lm}\label{lm PhiDoubleSpanTame}
The space $\varPhi'_d$
is spanned by the tame double $\phi\text{--words}$.
\end{lm}

\proof
The proof is in seven steps.
\begin{enumerate}
\item\label{item three}
By Lemma~\ref{lm IntersecT}, identity~\eqref{eq Doub-Trip Phi} gets in $\mathfrak M$ the form
$$
\linodd{n}{1}{2}{3}{\crc{4}}{\crc{5}}{\crc{6}}=0.
$$
Consequently,
in view of the cyclic property of~$\phi$,
a double
$\phi\text{--word}$
$\linodd{n}{1}{2}{3}{4}{5}{6}$
satisfies the assertion of Lemma~\ref{lm J-triples} for the variables
$$
x_{6}=x,\quad
x_{5}=y,\quad
x_{4}=z_1,\dots,\quad
x_{1}=z_4.
$$
Throughout the proof,
referring item~\ref{item three},
we apply  Lemma~\ref{lm J-triples}, combining it with Lemmas~\ref{lm cosreg} and~\ref{lm SvPhi} with no comments.
\item\label{item one}
Following the procedure of Lemma~\ref{lm PhiSpanTame}, one can prove that
every double $\phi\text{--word}$ can be represented as a linear combination of double
$\phi\text{--words}$ with the origins of tame $\phi\text{--words}$.
\item\label{item onee}
By item~\ref{item one} and Lemma~\ref{lm DoublePhiSpanNormal},
$\varPhi'_d$
is spanned by the normalized double $\phi\text{--words}$
with the origins of tame $\phi\text{--words}$.
\item\label{item two}
Let
$f\in\varPhi'_d$
be a wild normalized double $\phi\text{--word}$
with an origin of some tame
$\phi\text{--word}$.
Then $f$ has one of the following forms for $i<j<k$:
$$
\nlinodd{n}{1}{2}{d}{k}{j},\quad
\nlinodd{n}{1}{i}{k}{d}{j},\quad
\nlinodd{n}{1}{j}{k}{d}{i}.
$$
By item~\ref{item onee}, it suffices to prove that
$f$
can be expressed linearly with tame double $\phi\text{--words}$.
\item\label{item four}
If
$f=\nlinodd{n}{1}{2}{d}{k}{j}$ and
$j>m=\min\bigl(\{3,\dots,d-1\}\setminus \{k\}\bigr)$,
then by item~\ref{item three}, we have
$$
f=-\nlinodd{n}{1}{2}{d}{j}{k}-\linodd{n}{1}{2}{d}{m}{j}{k}.
$$
Here, the first summand is tame and the second one is either tame or, after normalization
by Lemma~\ref{lm DoublePhiSpanNormal}, gets the form
$$
\linodd{n}{1}{2}{d}{m}{j}{k}=
\pm \nlinodd{n}{1}{2}{d}{m}{j} \pm \nlinodd{n}{1}{2}{d}{m}{k},
$$
where the both summands are tame by the choice of $m$.
\item\label{item five}
If
$f=\nlinodd{n}{1}{i}{k}{d}{j}$, where
$j>m=\min\bigl(\{2,\dots,d-1\}\setminus \{i,k\}\bigr)$,
then by item~\ref{item three}, we obtain
$$
f=-\nlinodd{n}{1}{i}{j}{d}{k}-\linodd{n}{1}{i}{m}{d}{j}{k}.
$$
Again, the first summand is tame and if $i<m$, then the second summand is either tame or can be
normalized, by Lemma~\ref{lm DoublePhiSpanNormal}, up to two tame double $\phi\text{--words}$.
Otherwise, for $i>m$, using~\eqref{eq in-cyc}, we get
$$
\linodd{n}{1}{i}{m}{d}{j}{k}=-\linodd{n}{1}{m}{i}{d}{j}{k}-\linodd{n}{1}{m}{d}{i}{j}{k}.
$$
Now, as above, the normalization of the first summand gives two tame double $\phi\text{--words}$
and the second summand, by item~\ref{item one},
can be represented as a linear combination of double
$\phi\text{--words}$ with the origins of tame $\phi\text{--words}$ of type~2).
Consequently, by Lemma~\ref{lm DoublePhiSpanNormal} and item~\ref{item four}, these double
$\phi\text{--words}$ are also linear combinations of tame double
$\phi\text{--words}$.
\item\label{item six}
Finally, if
$f=\nlinodd{n}{1}{j}{k}{d}{i}$, where
$i>m=\min\bigl(\{2,\dots,d-1\}\setminus \{j,k\}\bigr)$,
then by item~\ref{item three}, we get
$$
f=-\nlinodd{n}{1}{i}{k}{d}{j}-\linodd{n}{1}{m}{k}{d}{i}{j}.
$$
The normalization of the second summand and the application of item~\ref{item five} complete the proof.
\hfill$\square$
\end{enumerate}

\subsection{Stable basis monomials and basis words}%
\label{SubSec:StableBasisMonomialsAndBasisWords}

A  basis monomial
$w\in\pdvv$
is called
\textit{stable}
if $w\notin{\mathcal P}^{(1,1)}_{d}(\M)$
or the minimal of the indices of the variables of its origin~$w_0$
is less then the minimal index in its formative operator~$\obr{w}$.

\textit{Basis words}
are all elements of
${\mathcal P}^{(1,1)}_{d}(\mathfrak M)$
of the following types:

1)\enskip $\psi\text{--words}$ of the form
$\psi_n\left(1^{\left<i\right>},2^{\left<i\right>},3^{\left<i\right>},4^{\left<i\right>},i\right)$,
$1\leqslant i\leqslant d$;

2)\enskip normalized tame $\phi\text{--words}$;

3)\enskip tame double $\phi\text{--words}$ for $d\geqslant7$.

\begin{lm}\label{lm PhiSpanBasis}
The space
${\mathcal P}^{(1,1)}_{d}(\mathfrak M)$
is spanned by its stable basis monomials and basis words.
\end{lm}

\proof
First we stress that
by the definition of $\psi\text{--word}$, taking into account Lemma~\ref{lm cosnach},
we can represent every nonstable basis monomial as a linear combination of
two stable basis monomials and a $\psi\text{--word}$ of the form
$f=\psi_n\left(1^{\left<i\right>},j,k,\ell,i\right)$.
Then combining Lemmas~\ref{lm cosreg} and~\ref{lm SvPsi} with identity~\eqref{eq cos-psi},
it is not hard to verify that
$f$
is skew-symmetric w.r.t. all its variables, except~$x_{1^{\left<i\right>}}$,
modulo
$\varPhi_{d}$
and linear combinations of stable basis monomials.
Consequently, every nonstable basis monomial can be expressed linearly with a stable one and
a basis word of type~1) modulo $\varPhi_{d}$.
Thus to complete the proof in the case $d=5$ it suffices to refer Lemma~\ref{lm PhiSpanTame}.
Further, for $d\geqslant7$, it follows from Lemmas~\ref{lm PhiSpanNormal} and~\ref{lm PhiSpanTame} that
$\varPhi_{d}$
is spanned by the basis words of type~2)
modulo
$\varPhi'_{d}$.
Finally, by Lemma~\ref{lm PhiDoubleSpanTame},
the basis words of type~3) are linear generators of
$\varPhi'_{d}$.
\endproof

\section{Proof of the Theorem}%
\label{Sec:ProofOfTheTheorem}

First we stress that the polynomial
$\Delta g_n$
of form~\eqref{eq Lin gn}
is regular and,
consequently, by Lemma~\ref{lm AddBasPd},
$\Delta g_n\neq 0$ in~$\Ag$.
Thus,
$\mathfrak M$
is a proper subvariety of~$\M$.
Furthermore, by Lemma~\ref{lm IntersecT}, we have
${\left(g_n\right)}^{\mathrm T}\cap{\left(g_m\right)}^{\mathrm T}=\{0\}$
for $n\neq m$.
Therefore, $\mathfrak M$ is non-finitely based.

\medskip

\begin{lm}\label{lm Gv2Vs}
The Grassmann algebra
$\Gv{\M}{2}$
lies in $\mathfrak M$.
\end{lm}

\proof
Consider the linearization
$f=\Delta g_n$
of form~\eqref{eq Lin gn}.
By Lemma~\ref{lm IntersecT}, it suffices to verify that
$f=0$
on every set of generators of~$\Gv{\M}{2}$.
Indeed, if we substitute the variables of
$\Tilde{f}$
by the generators of the free $\M\text{-superalgebra}$
$\svobs{\M}{u_1,u_2}$,
then at least two of the variables~$x_1,x_2,x_3$ get the same value.
Consequently, in view of the odd parity of
$u_1,u_2$,
the linear combination
$\Tilde{f}$
contains with every its monomial
$\alpha w$ $(\alpha=\pm1)$
the monomial
$-\alpha w$.
Hence,
$\Tilde{f}=0$
in
$\svobs{\M}{u_1,u_2}$.
\endproof

\medskip

Let $\A'=\A'_0\oplus\A'_1$ be a superalgebra
$$
\A'_0=\F\cdot a,\quad\;
\A'_1=\F\cdot v+\F\cdot w+\F\cdot y+\F\cdot z
$$
such that all nonzero products of its basis elements are the following:
$$
\t{z}{z}=a,\quad\t{y}{a}=v,\quad\t{z}{a}=w,\quad\t{y}{v}=\t{v}{y}=a.
$$
The direct verification shows that all the statements of
Sec.~\ref{Sec:VspomogatelnyeSuperalgebry} formulated for $\A$
hold for~$\A'$ as well.
Moreover,
$\A'$ satisfies all the Lemmas of Sec.~\ref{Sec:AddBasPdv}, except Lemma~\ref{lm raspad0},
that is true only for~$\varepsilon=0$.
Consequently,
taking into account Lemma~\ref{lm Gv2Vs},
we obtain that
$\A'$
is an
$\mathfrak M\text{-superalgebra}$
such that
all the stable basis monomials, not lying in
${\mathcal P}^{(1,1)}_{d}(\mathfrak M)$, for $d\geqslant5$,
are linearly independent on~$\Gob{\A'}$.
Therefore,
in view of  Lemmas~\ref{lm TAgMultlin} and~\ref{lm PhiSpanBasis},
to prove the Theorem it suffices to verify the linear independence of
the stable basis monomials and the basis words of
${\mathcal P}^{(1,1)}_{d}(\mathfrak M)$
on~$\Gob{\A'}$.

We set, as above,
$d=2n+3$, $n\in\mathbb N$
and also assume
$\varPhi'_5=0$.

\smallskip
\begin{lm}\label{lm Q-lin-nez}
If
$f=0$
in
$\Gob{\A'}$
for some
$f\in{\mathcal P}^{(1,1)}_{d}(\mathfrak M)$,
then
$f\in\varPhi_d$.
\end{lm}

\proof
By Lemma~\ref{lm PhiSpanBasis}, we can write down
$f$
in the form
$$
f\equiv\sum_{k=1}^{d}f_k\opL_{x_k}\pmod{\varPhi_d}
$$
such that
$$
f_k=\sum_{i=2}^{d-1}
\alpha_{i,k}
\left(x_1x_i\right){(LR)}^{n}
+\beta_{k}\,\cicl{1}{2}{3}{4}{(LR)}^{n-1},
$$
where
$\alpha_{i,k},\beta_{k}\in\F$.
By
$\Tilde{f}_{i,k}$
denote the value taken by the superpolynomial
$\Tilde{f}$
on the following elements of $\A'$:
$$
x_{i^{\left<k\right>}}=v,\quad
x_k=z,
\quad x_j=y
\enskip
\text{ for}\enskip
j\in\left[1,d\right] \setminus\left\{i^{\left<k\right>},k\right\}.
$$
While calculating the values
$\Tilde{f}_{i,k}$,
we take into account~\eqref{eq phi-kv}.
Suppose
$\alpha_{i,k}\neq 0$ for $i\geqslant4$;
then $\Tilde{f}_{i,k}$ is proportional to the element
$$
\left(\t{y}{v}\right){\left(L_yR_y\right)}^{n}L_{z}
=a{\left(L_yR_y\right)}^{n}L_{z}=\t{z}{a}=w\neq0.
$$
Otherwise, in view of Lemma~\ref{lm cosnach},
$f_k$
can be rewritten in the form
$$
f_k=\bigl(\alpha_{k}\left(x_1x_2\right)L_{x_3}
+\beta_{k}\left(x_2x_3\right)L_{x_1}
+\gamma_{k}\left(x_3x_1\right)L_{x_2}\bigr)R{(LR)}^{n-1},
$$
where
$\alpha_k=\beta_{k}+\alpha_{2,k}$
and
$\gamma_k=\beta_{k}-\alpha_{3,k}$.
If $\alpha_k,\beta_k,\gamma_k$
are not simultaneously zero, then
by virtue of the restriction
$\chr{\F}\neq 2$,
one of their pairwise sum is not zero as well.
Suppose, for example, that
$\alpha_k+\beta_k\neq0$;
then we have
$$
{\textstyle\frac{1}{\alpha_k+\beta_k}\Tilde{f}_{2,k}}
=\left(\alpha_k\,\t{y}{v}+\beta_k\,\t{v}{y}\right)
{\left(L_yR_y\right)}^{n}L_{z}
=a{\left(L_yR_y\right)}^{n}L_{z}=\t{z}{a}=w\neq0.
$$
Thus all the scalars in the considered expression of $f$ modulo $\varPhi_d$ prove to be zeros.
\endproof

\smallskip
\begin{lm}\label{lm Qs-lin-nez}
If
$f=0$
in
$\Gob{\A'}$
for some
$f\in\varPhi_d$,
then
$f\in\varPhi'_d$.
\end{lm}

\proof
By Lemmas~\ref{lm PhiSpanNormal} and~\ref{lm PhiSpanTame}, we can write down
$f$
in the form
$$
f\equiv
\sum_{j=3}^{d-1}
\sum_{i=2}^{j-1}
\alpha_{i,j}\,
\nlinod{n}{1}{i}{j}{d}
+\beta_{j}\,
\sum_{j=3}^{d-1}
\nlinod{n}{1}{2}{d}{j}
\pmod{\varPhi'_d},
$$
where
$\alpha_{i,j},\beta_j\in\F$.
By
$\Tilde{f}_{i,j}$
denote the value of
$\Tilde{f}$
on the following elements of $\A'$:
$$
x_i=x_j=z,\quad
x_k=y
\enskip
\text{ for}\enskip
k\in\left[1,d\right] \setminus\left\{i,j\right\}.
$$
Suppose
$\alpha_{i,j}\neq 0$;
then $\Tilde{f}_{i,k}$ is proportional to the element
$$
z^2{\left(L_yR_y\right)}^{n-1}L_y
=a{\left(L_yR_y\right)}^{n-1}L_y=\t{y}{a}=v\neq0.
$$
Otherwise, $f$ gets the form
$$
f\equiv
\sum_{j=3}^{d-1}
\beta_{j}\,
\nlinod{n}{1}{2}{d}{j}
\pmod{\varPhi'_d}.
$$
In the case $\beta_{j}\neq 0$, by the similar calculations, we obtain
$
\Tilde{f}_{d,j}\neq0.
$

Therefore, $f\in\varPhi'_d$.
\endproof

\medskip

\begin{lm}\label{lm IV-lin-nez}
If
$f=0$
in
$\Gob{\A'}$
for some
$f\in\varPhi'_d$,
then
$f=0$
in
$\mathfrak M$.
\end{lm}

\proof
By Lemma~\ref{lm PhiDoubleSpanTame}, we can write down
$f$
in the form
$$
f=
\sum_{k=4}^{d-1}
\sum_{j=3}^{k-1}
\sum_{i=2}^{j-1}
\alpha_{i,j,k}\,
\nlinodd{n}{1}{i}{j}{d}{k}
+
\sum_{k=4}^{d-1}
\sum_{j=3}^{k-1}
\beta_{j,k}\,
\nlinodd{n}{1}{2}{d}{j}{k},
$$
where
$\alpha_{i,j,k},\beta_{j,k}\in\F$.
By
$\Tilde{f}_{i,j,k}$
denote the value of
$\Tilde{f}$
on the following elements of $\A'$:
$$
x_i=x_j=x_k=z,\quad
x_\ell=y
\enskip
\text{ for}\enskip
\ell\in\left[1,d\right] \setminus\left\{i,j,k\right\}.
$$
Assume that
$\alpha_{i,j,k}\neq 0$;
then $\Tilde{f}_{i,j,k}$ is proportional to the element
$$
z^2{\left(L_yR_y\right)}^{n-1}L_z=\t{z}{a}=w\neq0.
$$
Otherwise, $f$ gets the form
$$
f=
\sum_{k=4}^{d-1}
\sum_{j=3}^{k-1}
\beta_{j,k}\,
\nlinodd{n}{1}{2}{d}{j}{k}.
$$
In the case $\beta_{j,k}\neq 0$, by the similar arguments, we obtain
$
\Tilde{f}_{d,j,k}\neq0.
$

Therefore all the scalars in $f$ are zeros.
\endproof

\medskip

Lemmas~\ref{lm Gv2Vs}--\ref{lm IV-lin-nez} yield\,
$
\mathfrak M=\vr\Gob{\A'}=\vr\Gv{\M}{2}.
$
Theorem is proved.
\medskip
\begin{rem}
{\upshape
The variety
$\mathfrak M$
doesn't satisfy the condition of minimality, i.~e. there are proper subvarieties of
$\mathfrak M$
that are also non-finitely based.
For instance, it follows from the proof of the Theorem that
by setting all the double
$\phi\text{--words}$
equal to zero we distinguish the proper non-finitely based subvariety of
$\mathfrak M$.}
\end{rem}

\subsection*{Acknowledgments}%
This article was carried out at the Department of Mathematics and Statistics (IME) of the University of S\~ao Paulo
as a part of the author's post-doc project supported by the FAPESP (2010/51880--2)
under the supervision of Prof. I.~P.~Shestakov.
The author is very grateful to his supervisor and to the IME for the kind hospitality and the creative atmosphere.
The author is also very thankful to Prof. S.~V.~Pchelintsev for his suggesting the problem
and for the useful discussions on the obtained results.

\end{document}